\documentclass[10pt, english,preprint,jcp,superscriptaddress]{revtex4-1}
\usepackage[T1]{fontenc}
\usepackage[latin9]{inputenc}
\usepackage{xcolor}
\usepackage{pdfcolmk}
\usepackage{float}
\usepackage{subfig}
\usepackage{amsmath}
\usepackage{amssymb}
\usepackage{graphicx}
\usepackage{hyperref}

\newcommand{\ud}{\,\mathrm{d}}

\newcommand{\RR}{\mathbb{R}}

\newcommand{\Or}{\mathcal{O}}

\newcommand{\bd}[1]{\boldsymbol{#1}}
\newcommand{\wt}[1]{\widetilde{#1}}

\IfFileExists{mathabx.sty}%
  {\DeclareFontFamily{U}{mathx}{\hyphenchar\font45}%
   \DeclareFontShape{U}{mathx}{m}{n}{<->mathx10}{}%
   \DeclareSymbolFont{mathx}{U}{mathx}{m}{n}%
   \DeclareFontSubstitution{U}{mathx}{m}{n}%
   \DeclareMathAccent{\widebar}{0}{mathx}{"73}%
}{%
  \PackageWarning{mathabx}{%
    Package mathabx not available, therefore\MessageBreak substituting
    widebar with overline\MessageBreak }%
  \newcommand{\widebar}[1]{\overline{#1}}%
}

\newcommand{\eps}{\epsilon}

\newcommand{\abs}[1]{\lvert#1\rvert}

\newcommand{\norm}[1]{\left\lVert#1\right\rVert}

\newcommand{\braket}[2]{(#1\vert#2)}

\newcommand{\aux}{\mathrm{aux}}

\begin{document}

\title{Compression of the electron repulsion integral tensor in tensor
  hypercontraction format with cubic scaling cost}

\author{Jianfeng Lu} \affiliation{Departments of Mathematics, Physics,
  and Chemistry, Duke University, Box 90320, Durham, NC 27708}
\email{jianfeng@math.duke.edu}

\author{Lexing Ying} \affiliation{Department of Mathematics and
  Institute for Computational and Mathematical Engineering, Stanford
  University, 450 Serra Mall, Bldg 380, Stanford, CA 94305}
\email{lexing@stanford.edu}

\date{\today}

\begin{abstract}
  Electron repulsion integral tensor has ubiquitous applications in
  electronic structure computations. In this work, we propose an
  algorithm which compresses the electron repulsion tensor into the
  tensor hypercontraction format with $\Or(nN^2\log N)$ computational
  cost, where $N$ is the number of orbital functions and $n$ is the
  number of spatial grid points that the discretization of each
  orbital function has. The algorithm is based on a novel strategy of
  density fitting using a selection of a subset of spatial grid points
  to approximate the pair products of orbital functions on the whole
  domain.
\end{abstract}

\maketitle

\section{Introduction}

Given a set of orbital functions $\{\psi_i\} \subset L^2(\RR^3)$, the
four-center two-electron repulsion integrals
\begin{equation}
  \braket{ij}{kl} = \iint_{\RR^3\times \RR^3} \frac{\psi_i(x) \psi_j(x) \psi_k(y) \psi_l(y)}{\abs{x - y}} \ud x \ud y
\end{equation}
are universally used in many electronic structure theories, such as
Hartree-Fock, density functional theory (DFT), RPA, MP2, CCSD, and GW. As a
result, a key step to accelerate ab initio computations in quantum
chemistry and materials science is to get an efficient representation
of the electron repulsion integral tensor.

One of the most popular methods for compressing the electron repulsion
integral is the density fitting approximation. This method, also known
as resolution of identity approach \cite{Ren_etal:12,
  DunlapConnollySabin:79, Schutz_etal:10, SodtSubotnikHeadGordon:06,
  Vahtras:93, Weigend:98}, amounts to representing pair products of
orbital functions $\psi_i(x) \psi_j(x)$ in terms of a set of auxiliary
basis functions
\begin{equation}\label{eq:densityfitting}
  \rho_{ij}(x) := \psi_i(x) \psi_j(x) \approx \wt{\rho}_{ij}(x) = \sum_{\mu} C_{ij}^{\mu} P_{\mu}(x),
\end{equation}
where $\mu = 1, 2, \ldots, N_{\aux}$ labels the auxiliary basis
functions. The auxiliary basis functions are constructed either
explicitly (e.g., a set of Gaussian-type atom-centered basis
functions) or implicitly by using singular value decomposition on the
overlap matrix of the set of $N^2$ functions $\rho_{ij}(x)$
\cite{Foerster:08, Foerster:11}.

After the auxiliary basis functions are determined, a least square
fitting is used to determine the coefficient $C_{ij}^{\mu}$. 
When the standard $L^2$ metric is used in the least
square fitting, one obtains 
\begin{align}
  & C_{ij}^{\mu} = \sum_{\nu} \langle ij \vert \nu \rangle S_{\nu\mu}^{-1}, \\
  & \braket{ij}{kl} \approx \sum_{\mu\mu'\nu\nu'} \langle ij \vert \nu \rangle S_{\nu \mu}^{-1}  V_{\mu\mu'} S_{\mu'\nu'}^{-1} \langle \nu' \vert kl \rangle
\end{align}
with the short-hand notations
\begin{align}
  & \langle ij \vert \nu \rangle = \int \psi_i(x) \psi_j(x) P_{\nu}(x) \ud x,  \\
  & S_{\mu\nu} = \int P_{\mu}(x) P_{\nu}(x) \ud x, \qquad \text{and}\\
  & V_{\mu\nu} = \iint \frac{P_{\mu}(x) P_{\nu}(y)}{\abs{x - y}}
  \ud x \ud y. 
\end{align}
It is also possible to use the Coulomb weight in the least square fitting, which leads to 
\begin{align}
  & C_{ij}^{\mu} = \sum_{\nu} \braket{ij}{\nu} V_{\nu\mu}^{-1}, \\
  & \braket{ij}{kl} \approx \sum_{\mu\nu} \braket{ij}{\mu} V_{\mu\nu}^{-1} 
  \braket{\nu}{kl}
\end{align}
with the short-hand notation
\begin{equation}
   \braket{ij}{\nu} = \iint \frac{\psi_i(x) \psi_j(x) P_{\nu}(y)}{\abs{x - y}} \ud x \ud y.
\end{equation}
A closely related idea to density fitting is the incomplete Cholesky
decomposition of the electron repulsion integrals
\cite{BeebeLinderberg:77, Koch:03}:
\begin{equation}
  \braket{ij}{kl} \approx \sum_{\mu=1}^{M} L_{ij}^{\mu} L_{kl}^{\mu},
\end{equation}
where $L_{ij}^{\mu}$ are numerically obtained Cholesky vectors.  The
cost of getting the resolution of identity approximation, assuming
$\Or(N)$ auxiliary basis functions, is $\Or(N^4)$, where $N$ is the
number of orbital functions. Other methods for the electron repulsion
integral tensor include multipole moment approaches
\cite{GreengardRokhlin:87, White:94, White:96, Strain:96} and
pseudospectral representation \cite{Friesner:85, Martinez:92,
  MartinezCarter:95}.



More recently, the tensor hypercontraction of the electron repulsion
integral have been proposed in \cite{HohensteinParrishMartinez:12,
  ParrishHohensteinMartinez:12, ParrishHohenstein:13}, which aims at an approximation of the
electron repulsion integral tensor as
\begin{equation}
  \braket{ij}{kl} \approx \sum_{\mu\nu} X_i^{\mu} X_j^{\mu} Z^{\mu\nu}
  X_k^{\nu} X_l^{\nu},
\end{equation}
where $\mu, \nu$ are the indices for the decomposition. The factor $X$
is taken to be the weighted collocation matrix arises from numerical
quadrature of the electron repulsion integral and $Z$ is determined by
a least square procedure. The computational cost of obtaining the
approximation is either $\Or(N^5)$ when direct quadrature of electron
repulsion integral is used or $\Or(N^4)$ with the help of density
fitting procedure.  The tensor hypercontraction opens doors to
efficient algorithms for several electronic structure theories, see
e.g., \cite{Hohenstein:13, HohensteinParrishMartinez:12,
  ParrishHohensteinMartinez:12, Parrish:14, Shenvi:13, Shenvi:14}.


In this work, we propose an $\Or(n N^2 \log N)$ algorithm to get the
tensor hypercontraction of the electron repulsion integral. It is
based on an approximation of $\rho_{ij}(x)$ similar to
\eqref{eq:densityfitting}, but with the key advantage that the
coefficient $C_{ij}^{\mu}$ has separate dependence on the indices $i$
and $j$. Such an approximation is achieved by an interpolative
decomposition which chooses selected grid points $x_{\mu}$ to
interpolate the pair product density $\rho_{ij}$. This is different
from the usual density fitting strategy with a predetermined set of
auxiliary basis functions. In this sense, our algorithm tries to find
an optimal set of the auxiliary basis functions, such that the tensor
hypercontraction format can be immediately obtained.

\section{Algorithm}

Our algorithm is based on the randomized column selection method for
low-rank matrix, recently developed in \cite{Liberty:07,
  Woolfe:08}. For an $m \times n$ matrix $A$, the column selection
method looks for an interpolative decomposition to approximate $A
\approx C P$ such that the discrepancy $\norm{A - CP}$ is minimized,
where $C$ is an $m \times c$ matrix consists of $c$ columns of $A$ and
$P$ is a $c \times n$ matrix. The interpolative decomposition based on
randomized column selection has recently been used for finding Wannier
functions given a set of eigenfunctions in Kohn-Sham density
functional theory \cite{DamleLinYing} by one of the authors. Here we
demonstrate the power of the interpolative decomposition in the
context of compressing electron repulsion integral tensor.

In our context, we will apply the column selection method on
$\rho_{ij}(x)$ which is viewed as an $(N^2) \times n$ matrix, where
$N$ is the number of orbitals $\psi_i$ and $n$ is the total number of
spatial grid points, \textit{i.e.}, we will view the pair $(ij)$ as
the row index and the grid point $x$ as the column index of the
matrix. We remark that while we will treat $x$ as a spatial grid
throughout the presentation for definiteness, in other words, we have
assumed a real space discretization of $\psi_i$, in fact, it is also
possible to extend the algorithm to other discretizations, e.g.,
atomic orbitals, by using the idea proposed in pseudospectral
representation \cite{Friesner:85, Martinez:92, MartinezCarter:95}. Let
us emphasize that the choice of the spatial quadrature grid $x$ is
completely general in our methods.

The column selection then amounts to choose a number of
spatial grid points, denoted as $x_{\mu}$, $\mu = 1, \ldots,
N_{\aux}$, such that $\rho_{ij}(x)$ is approximated as
\begin{equation}\label{eq:columnselection}
  \rho_{ij}(x) \approx \sum_{\mu} \rho_{ij}(x_{\mu}) P_{\mu}(x) = \sum_{\mu} \psi_i(x_{\mu}) \psi_j(x_{\mu}) P_{\mu}(x). 
\end{equation}
This should be compared with the approximation in the density fitting
\eqref{eq:densityfitting}: Here $\psi_i(x_{\mu}) \psi_j(x_{\mu})$
plays the role of the coefficient $C_{ij}^{\mu}$ in
\eqref{eq:densityfitting}, which is the key feature of the
interpolative decomposition approximation. To avoid possible
confusion, unlike what is commonly involved in conventional density
fitting approaches, the approximation \eqref{eq:columnselection} is
not a quadrature formula, it should be understood as an
interpolation. In particular, this should be distinguished from the
flavor of tensor hypercontraction known as X-THC in
\cite{ParrishHohenstein:13}, which is essentially a Gaussian
quadrature formula for the overlap integrals.

The approximation \eqref{eq:columnselection}  has a clear advantage
that the dependence on $i$ and $j$ are separated as a result of using
the selected columns to approximate the whole matrix. Indeed, assuming
such an approximation \eqref{eq:columnselection}, the electron
repulsion integral tensor then becomes
\begin{equation}
  \braket{ij}{kl} \approx \sum_{\mu\nu} \psi_i(x_{\mu}) \psi_j(x_{\mu}) 
  V_{\mu\nu} \psi_k(x_{\nu}) \psi_j(x_{\nu}). 
\end{equation}
Hence, we immediately arrive at the tensor hypercontraction format of
the electron repulsion integral tensor without further approximation!
The only extra step is to calculate $(\mu\vert\nu)$, which can be done
efficiently using fast Fourier transform (FFT). 

It remains to show how an approximation as \eqref{eq:columnselection}
can be efficiently obtained. As opposed to the density fitting
approach, the central focus in our algorithm is the selection of
columns. After $N_{\aux}$ grid points $x_{\mu}$ are determined, the
auxiliary basis functions $P_{\mu}$ follow from least squares
fitting. To find the suitable subset of columns, a pivoted QR
algorithm \cite{GolubVanLoan} is used on a random projection of
$\rho_{ij}(x)$. In more details, the algorithm for the column
selection consists of the following steps, given $\rho_{ij}(x)$ and an
error threshold $\eps$.
\begin{enumerate}
\item Reshape $\rho_{ij}(x)$ into an $(N^2) \times n$ matrix
  by combining $(ij)$ as a single index: 
  \begin{equation}
    \varrho_{(i-1)N + j}(x) = \rho_{ij}(x), 
  \end{equation}
  where the index of $\varrho$, which will be denoted as $I$ in the
  following, goes from $1$ to $N^2$;
\item Random Fourier projection of $\varrho_I(x)$:
  \begin{enumerate}
  \item Compute for $\xi = 1, \ldots, N^2$ the discrete Fourier
    transform
    \begin{equation}
      \mathfrak{M}_{\xi}(x) = \sum_{I=1}^{N^2} e^{-2 \pi \sqrt{-1} I \xi / N^2} \eta_I \varrho_I(x), 
    \end{equation}
    where $\eta_I$ is a random unit complex number for each $I$.
  \item Choose a submatrix $M$ of $N^2 \times n$ matrix $\mathfrak{M}$
    by randomly choosing $r N$ rows. In practice, $r = 20$ is used in
    our implementation.
  \end{enumerate}
\item Compute the pivoted QR decomposition of the $rN \times n$ matrix
  $M$: $M E = QR$, where $E$ is an $n \times n$ permutation matrix, $Q$
  is a $rN \times rN$ unitary matrix, and $R$ is a $rN \times n$ upper
  triangular matrix with diagonal entries in decreasing order.

  \noindent Note that $ME$ amounts to a permutation of the columns of
  $M$.
\item Determine the number of auxiliary basis functions $N_{\aux}$,
  such that $ \abs{R_{N_{\aux}+1, N_{\aux}+1}} < \eps \abs{R_{1, 1}}
  \leq \abs{R_{N_{\aux}, N_{\aux}}}$, i.e., this is a thresholding of
  the diagonals of $R$ to the relative error threshold $\eps$.
\item Choose $x_{\mu}$, $\mu = 1, \ldots, N_{\aux}$ such that the
  $x_\mu$-column of $M$ corresponds to one of the first $N_{\aux}$
  columns of $M E$.
\item Denote $R_{1:N_{\aux}, 1:N_{\aux}}$ the submatrix of $R$ consists
  of its first $N_{\aux} \times N_{\aux}$ entries, and $R_{1:N_{\aux},
    \bd{:}}$ the submatrix consists of the first $N_{\aux}$ rows of $R$.
  Compute
  \begin{equation*}
    P = R_{1:N_{\aux}, 1:N_{\aux}}^{-1} R_{1:N_{\aux}, \bd{:}} \, E^{-1}.
  \end{equation*}
  Then each row of the $N_{\aux} \times n$ matrix $P$ gives an
  auxiliary basis function $P_{\mu}(x)$ for $\mu = 1, \ldots,
  N_{\aux}$.
\end{enumerate}

The computationally expensive steps of the above algorithm are Steps
2, 3, and 6. Step 2 takes $n$ times FFT of $N^2$ length vectors, and
hence has complexity $\Or(n N^2 \log N)$. Step 3 computes QR
decomposition of $M$, which has complexity $\Or(n N^2 )$. Step 6
involves inversion of an $N_{\aux} \times N_{\aux}$ matrix and multiply
the inverse with an $N_{\aux} \times n$ matrix, which has complexity
$\Or(N_{\aux}^3 + n N_{\aux}^2)$. Hence, the overall complexity of the
column selection is $\Or(n N^2 \log N)$, as $N_{\aux} = \Or(N)$. The
memory storage cost of the intermediate results is also $\Or(n N^2)$,
which is the same as the cost of storing each entry of $\rho_{ij}(x)$.

Note that the Fourier transform in Step 2 of the algorithm acts on the
index of the pair densities, but not the spatial grids. The Fourier
transform is used for the random projection. We emphasize again that
our algorithm does not rely on any particular choice of the spatial
grids.

\section{Numerical results}

Given a set of orbital functions $\{\psi_i\}$, we denote
$\wt{\rho}_{ij}$ the result of the approximation based on the column
selection method in the previous section. We measure the error in
two ways by using the $L^2$ metric and the Coulomb metric:
\begin{align}
  & e_{ij}^{(2)} = \biggl( \int \abs{\rho_{ij}(x) - \wt{\rho}_{ij}(x)}^2 \ud x \biggr)^2; \\
  & e_{ij}^{(c)} = \biggl( \iint \frac{(\rho_{ij} - \wt{\rho}_{ij})(x)
    (\rho_{ij} - \wt{\rho}_{ij})(y)}{\abs{x - y}} \ud x \ud y\biggr)^{1/2}. 
\end{align}
Note that the approximation error of the electron repulsion tensor can be controlled by $\max e_{ij}^{(c)}$, since we have 
\begin{equation}
  \begin{aligned}
    & \braket{ij}{kl} - \sum_{\mu\nu}  \psi_i(x_{\mu}) \psi_j(x_{\mu}) \braket{\mu}{\nu} \psi_k(x_{\nu}) \psi_l(x_{\nu}) \\
    & \quad = \iint \frac{\rho_{ij}(x) \rho_{kl}(y) -
      \wt{\rho}_{ij}(x)
      \wt{\rho}_{kl}(y)}{\abs{ x- y}} \ud x \ud y \\
    & \quad = \iint \frac{\rho_{ij}(x) (\rho_{kl}(y) - \wt{\rho}_{kl}(y))}{\abs{ x- y}} \ud x \ud y  \\
    & \quad \qquad + \iint \frac{( \rho_{ij}(x) - \wt{\rho}_{ij}(x))
      \wt{\rho}_{kl}(y)}{\abs{ x- y}} \ud x \ud y \\
    & \quad \leq \norm{\rho_{ij}}_C e_{kl}^{(c)} + e_{ij}^{(c)}
    \norm{\wt{\rho}_{kl}}_C, 
  \end{aligned}
\end{equation}
where the last inequality follows from the Cauchy-Schwartz inequality and $\norm{\cdot}_C$ stands for the Coulomb norm:
\begin{equation}
  \norm{f}_C = \biggl( \iint \frac{f(x) f(y)}{\abs{x - y}} \ud x \ud y\biggr)^{1/2}. 
\end{equation}

We first test the performance of the algorithm for an $1D$ toy problem
where the orbital functions are chosen to be the first $N$
eigenfunctions of a Hamiltonian operator $ H = - \tfrac{1}{2} \Delta +
V$, discretized on an interval rescaled to $[0, 1]$ with $n = 1024$
grid points, the periodic boundary conditions are used. To be consist
with the periodic boundary condition, we replace the bare Coulomb
interaction with the periodic Coulomb interaction to take into account
the interaction with periodic images. Taking $V$ to be a potential
randomly generated that consists of the first $128$ Fourier modes on
$[0, 1]$, we first diagonalize the discretized Hamiltonian to obtain
$\{\psi_i\}$ and then apply the column selection method. We test the
performance using different values of the threshold $\eps$ in Step 4
of the algorithm. The result is shown in Table~\ref{tab:1Derror},
where the dimensionless relative errors are defined to be
\begin{align}
  & \text{rel.~$2$-error} = \text{mean}(e_{ij}^{(2)}) / \text{mean} \norm{\rho_{ij}}_2; \\
  & \text{rel.~c-error} = \text{mean}(e_{ij}^{(c)}) / \text{mean} \norm{\rho_{ij}}_C,
\end{align}
where the average is taken with respect to the $N^2$ indices $(ij)$.
We observe that the error measured in both the $L^2$ metric and the
Coulomb metric is well controlled by the threshold $\eps$ with a small
number of auxiliary functions. Note that we have $N^2 = 16384$ pair of
orbitals in this example, while $10^{-5}$ relative error is achieved
with $N_{\aux}$ around $300$. We also note that the number of
auxiliary functions only increase mildly as we reduce the error
threshold.
\begin{table}[ht]
  \centering
  \begin{tabular}{|c|c|c|c|c|c|}
    \hline
    $\eps$ & $N_{\aux}$ & $\max e_{ij}^{(2)}$ & $\max e_{ij}^{(c)}$ & rel.~$2$-error & rel.~c-error  \\
    \hline
    1E-5 & 300 & 1.477E-7 & 9.154E-6 & 6.806E-6 & 1.051E-5 \\
    1E-6 & 324 & 1.095E-8 & 8.626E-7 & 9.747E-7 & 1.366E-6 \\
    1E-7 & 353 & 1.877E-9 & 2.035E-7 & 1.086E-7 & 1.610E-7 \\
    \hline
  \end{tabular}
  \caption{Error of the density fitting by column selection in $1D$ with $N = 128$ and $n = 1024$.}\label{tab:1Derror}
\end{table}

To test the computational complexity of the algorithm, we use a range
of $N$ and $n$ while keeping the same error threshold $\eps =
10^{-5}$. The timing results are shown in Table~\ref{tab:1Dtime}
together with the error of the fitting. The algorithm is implemented
using Matlab and the test is done on a single core on Intel Xeon CPU
X5690 \@ 3.47GHz. The timing matches very well with the complexity
$\Or(nN^2 \log N)$. The linear dependence of $N_{\aux}$ on $N$ is also
apparent. 
\begin{table}[ht]
  \centering
  \begin{tabular}{|c|c|c|c|c|c|}
    \hline
    $N$ & $n$ & $N_{\aux}$ & rel.~$2$-error & rel.~c-error & time \\
    \hline
    64 & 512 & 154 & 7.101E-6 & 1.534E-5 & 0.077s \\
    128 & 512 & 287 & 5.591E-6 & 3.472E-6 & 0.217s \\
    128 & 1024 & 304 & 1.011E-5 & 2.707E-5 & 0.467s \\
    256 & 1024 & 584 & 7.214E-6 & 6.268E-6 & 1.550s \\
    256 & 2048 & 593 & 1.089E-5 & 2.555E-5 & 4.244s \\
    512 & 2048 & 1156 & 5.355E-6 & 4.533E-6 & 17.881s \\
    \hline
  \end{tabular}
  \caption{Error and timing of the density fitting by column selection in $1D$ for fixed $\eps = 10^{-5}$ and different $N$ and $n$.}\label{tab:1Dtime}
\end{table}

To further test the algorithm in $3D$, we perform a $3D$
generalization of the numerical test in $1D$, where the orbital
functions are taken to be collection of eigenfunctions of a given
Hamiltonian operator. Here we take $n = 4096$ degrees of freedom for
each orbital $\psi_i$ and vary $N$ the number of orbital functions and
$\eps$ the error threshold to evaluate the performance of the
algorithm. The results are shown in Table~\ref{tab:3Derror}. We
observe that the relative error is still well controlled by the error
threshold $\eps$, while in $3D$ we need more auxiliary basis functions
compared to $1D$ case. For different $N$ and fixed error threshold
$\eps$, the number of auxiliary basis functions grows roughly linearly
with respect to $N$, confirming the scaling $N_{\aux} = \Or(N)$. Note
that in all cases, $N_{\aux}$ is much smaller compared with the total
possible pair of orbitals $N^2$. The computational time also agrees
well with $\Or(n N^2 \log N)$ (note that $n$ is fixed in this
example).
\begin{table}[ht]
  \centering
  \begin{tabular}{|c|c|c|c|c|c|c|}
    \hline
    $\eps$ & $N$ & $N_{\aux}$ & $\max e_{ij}^{(c)}$ & rel.~$2$-error & rel.~c-error & time \\
    \hline
    1E-4 & 32 & 303 & 1.863E-4 & 1.115E-4 & 6.290E-5 & 1.856s \\
    1E-5 & 32 & 358 & 2.462E-5 & 1.676E-5 & 8.966E-6 & 1.892s \\
    1E-6 & 32 & 407 & 3.429E-6 & 2.442E-6 & 1.280E-6 & 1.906s \\
    \hline
    1E-4 & 64 & 647  & 2.027E-4 & 1.126E-4 & 5.879E-5 & 6.403s \\
    1E-5 & 64 & 767  & 2.257E-5 & 1.648E-5 & 8.330E-6 & 6.429s \\
    1E-6 & 64 & 891  & 2.891E-6 & 2.484E-6 & 1.209E-6 & 6.431s \\ 
    \hline
    1E-4 & 128 & 1323 & 2.006E-4 & 1.018E-4 & 5.259E-5 & 20.046s \\
    1E-5 & 128 & 1516 & 1.913E-5 & 1.285E-5 & 6.367E-6 & 20.212s \\
    1E-6 & 128 & 1731 & 2.538E-6 & 1.661E-6 & 7.868E-7 & 20.497s \\
    \hline
  \end{tabular}
  \caption{Error and timing of the density fitting by column selection in $3D$ with 
    $n = 4096$.}\label{tab:3Derror}
\end{table}

\begin{figure}[ht]
  \centering 
  \subfloat{
    \includegraphics[width = 0.48\textwidth]{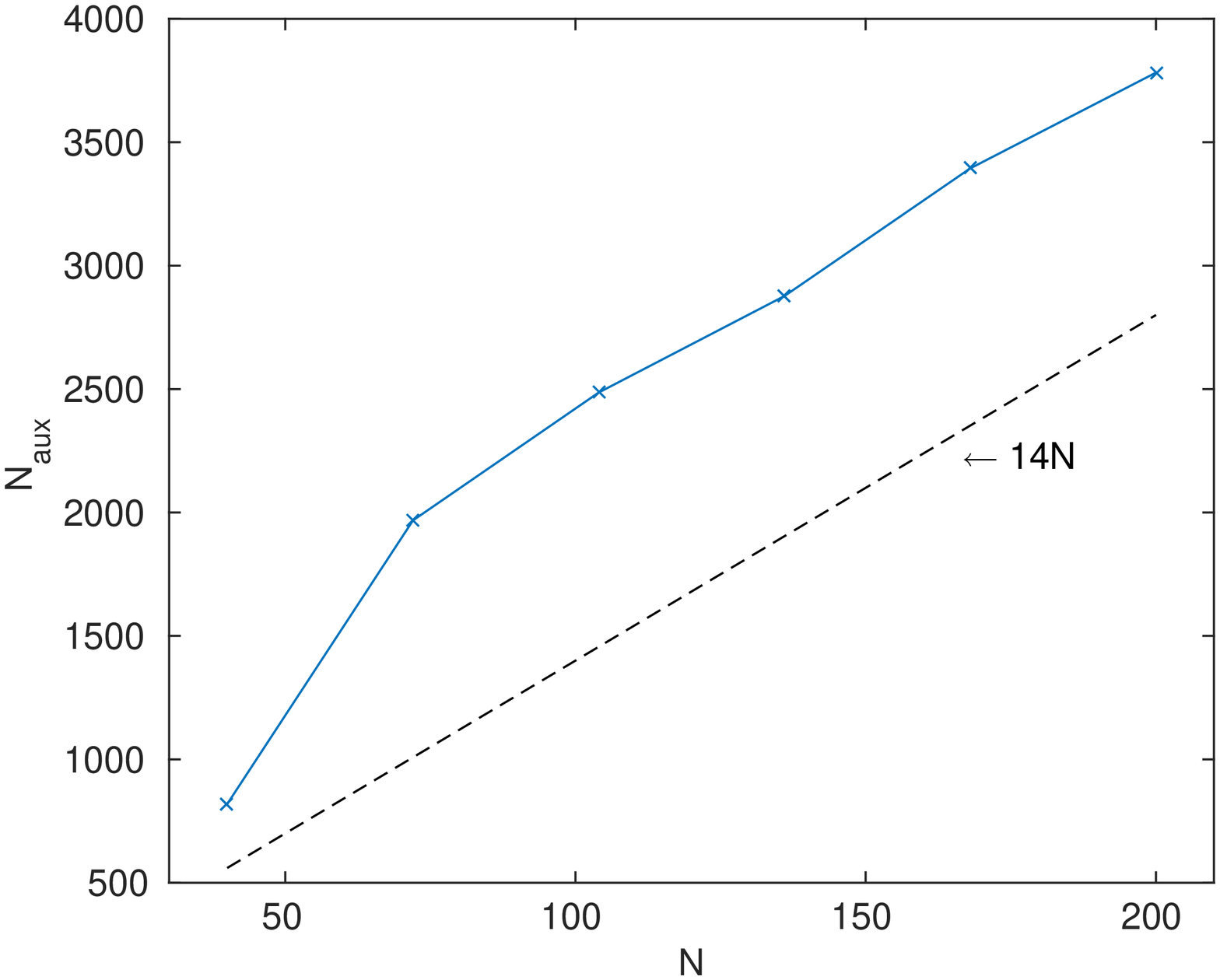}} \quad
  \subfloat{
    \includegraphics[width = 0.48\textwidth]{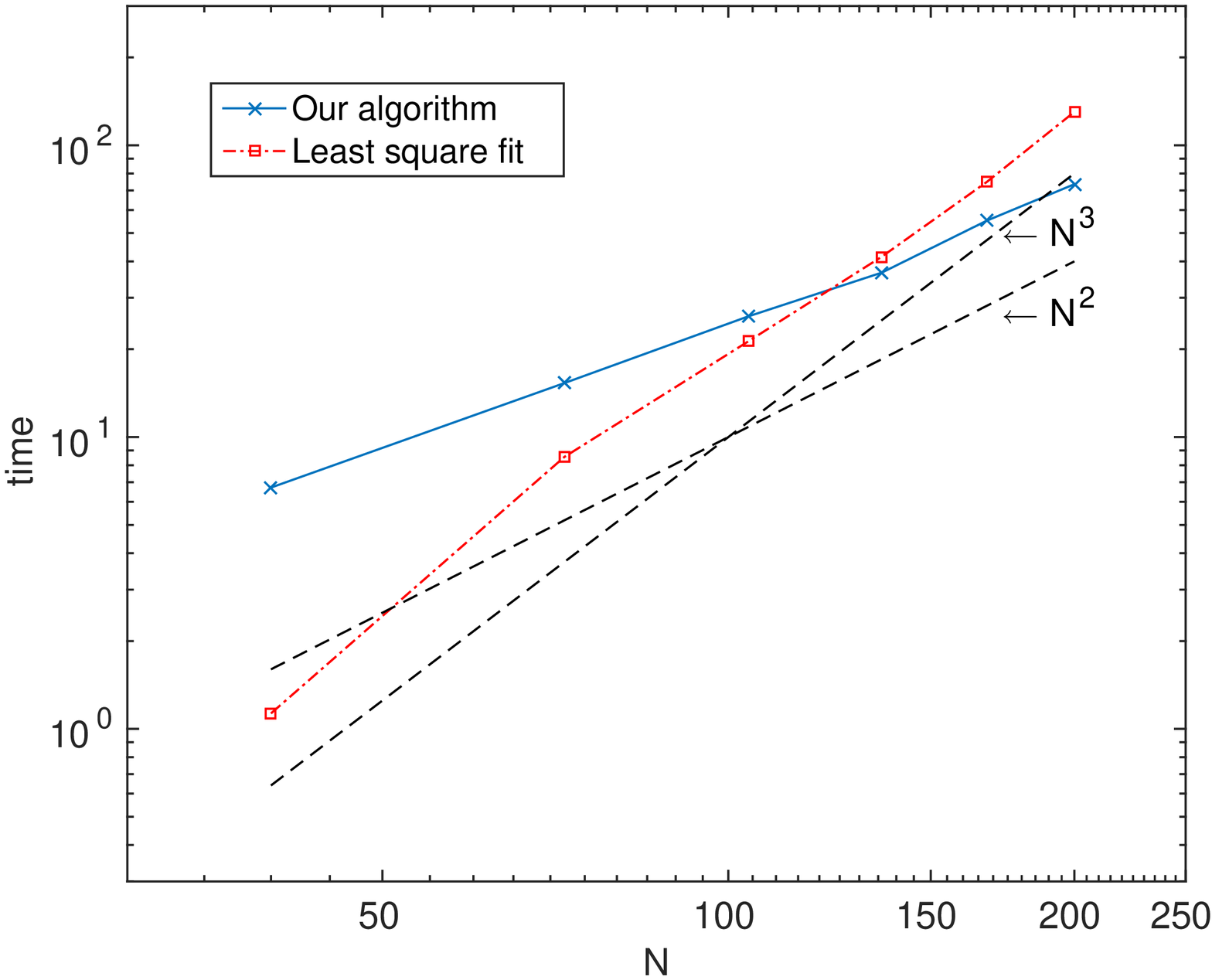}}
  \caption{(left) $N_{\aux}$ as a function of $N$. The dash line
    provides a reference of the linear dependence on $N$. (right)
    Comparison of timing of the current algorithm and the conventional
    density fitting based on least square fitting. The dash lines
    provide reference of quadratic and cubic dependence on $N$. \label{fig:time}}
\end{figure}

Finally, we consider a more realistic example based on the
implementation of the proposed algorithm in \textsf{KSSOLV}
\cite{kssolv}, a \textsf{MATLAB} toolbox for solving the Kohn-Sham
equations. For the test example, we choose two unit cells of a
graphene sheet (and hence consisting of $4$ carbon atoms) with
periodic boundary condition. Planewave is used for spatial
discretization with a fixed energy cutoff, and hence $n$ is fixed. We
take the first $N$ orbitals of the self-consistent Hamiltonian to be
the collection of orbitals for density fitting. The error threshold
$\eps$ is fixed for different $N$. The results are shown in
Table~\ref{tab:time}. To compare the algorithm with conventional
density fitting based on least square fitting with $L^2$ metric, we
also include the timing of the conventional calculation of the
coefficient based on the same auxiliary basis obtained in the proposed
algorithm. The comparison of timing is further illustrated in
Figure~\ref{fig:time} (right). Since $n$ is fixed in this test, our
algorithm scales as $N^2 \log N$ and the conventional density fitting
scales as $N^3$, which are clearly seen on the figure. Hence, for
large $N$, the current algorithm has lower computational cost, even
compared to density fitting, which is a preliminary step to get
hypercontraction format. The Figure~\ref{fig:time} (left) verifies the
linear scaling dependence of $N_{\aux}$ on $N$. We note that except
for a pre-asymptotic regime for small $N$, the linear dependence is
clear.

\begin{table}[ht]
  \centering
  \begin{tabular}{|c|c|c|c|c|c|}
    \hline
    $N$ & $N_{\aux}$ & rel.~$2$-error & rel.~c-error & time (proposed alg.) & time (least sq. fit)\\
    \hline
    8 & 36 & 5.658E-12 & 5.555e-12 & 1.148s & 0.0161s \\
    40 & 819 & 1.346E-4 & 6.571E-5 & 6.713s & 1.130s  \\
    72 & 1968 & 6.803E-4 & 2.922e-4  & 15.310s & 8.548s \\
    104 & 2486 & 3.939E-4 & 1.770E-4 & 25.890s & 21.239s \\
    136 & 2877 & 2.360E-4 & 1.119E-4 & 36.607s & 41.514s \\
    168 & 3394 & 8.068E-5 & 3.796E-5 & 55.244s & 75.074s \\
    200 & 3782 & 4.685E-5 & 2.163E-5 & 73.514s & 130.041s \\
    \hline
  \end{tabular}
  \caption{Error and timing of the density fitting by column selection (implemented in \textsf{KSSOLV}) for fixed $\eps = 10^{-5}$, $n = 10602$ and different $N$.}\label{tab:time}
\end{table}

\section{Discussion and conclusion}

The proposed cubic scaling algorithm for tensor hypercontraction
format of electron repulsion integral tensor is easy to implement and
can be easily incorporated into existing electronic structure
packages. Relatively small scale numerical tests are done in this
manuscript to demonstrate the effectiveness of the
algorithm. Applications to large scale electronic structure
calculations are the natural next steps.

The algorithm applies to general collection of orbital functions. In
particular, we do not assume any locality of the functions
$\{\psi_i\}$ in the algorithm. If a set of localized orbitals / basis
functions are considered, it is then possible to utilize the locality
to further reduce the computational cost. For instance, for a sub-collection of the orbitals, we may localize the  column selection to the support of them. 
It would then even possible
to reduce the computational scaling to $\Or(N)$ with controllable
error.  This is an important future direction that we plan to pursue.

We also remark that for the simplicity of the presentation, here we
have assumed that the orbital function $\psi_i$s are already
represented on a real space grid. We emphasize that the choice of the
spatial grid can be quite flexible. For example, if atomic orbital
discretization is used, one can first get a real space representation
using quadrature grids and then apply our algorithm. The computational
complexity depends on $n$, the number of spatial grid points, which in
practice will be a constant factor of $N$, while this prefactor might
be large. It would be interesting to explore algorithms that can work
directly with atomic orbital functions without first going to the real
space representation.

Also related to the previous point of changing basis functions. The
column selection method is designed with the discrepancy given by the
Frobenius norm, i.e., $L^2$ metric. While our numerical tests have
shown that the performance measured in error in either $L^2$ metric or
Coulomb metric is satisfactory, one observes that the error in Coulomb
metric is slightly larger than in the $L^2$ metric. It is therefore
interesting to ask whether the column selection can be done in Coulomb
metric directly. The natural idea of working on the Fourier domain
does not work, as the Fourier transform in $x$ will destroy the
separability of the dependence of the coefficients $C_{ij}^{\mu}$ on
which the algorithm crucially depends on. To avoid possible confusion,
let us emphasize that while the column selection uses $L^2$ metric,
the density fitting proposed by the current algorithm is actually
quite different from the RI-SVS density fitting (see e.g, the review
article \cite{Ren_etal:12}).

Finally, it would be interesting to explore fast algorithms for
quantum chemistry calculations based on the $\Or(nN^2\log N)$ algorithm for
tensor hypercontraction proposed here.

\medskip
\noindent\textbf{Acknowledgment.} J.L. would like to thank Weitao Yang
for helpful discussions. The work of J.L. is supported in part by the
Alfred P.~Sloan Foundation and the National Science Foundation under
grant DMS-1312659. The work of L.Y. is partially supported by the
National Science Foundation under grant DMS-0846501 and the
U.S. Department of Energy's Advanced Scientific Computing Research
program under grant DE-FC02-13ER26134/DE-SC0009409.

\bibliographystyle{plainnat}
\bibliography{dft}

\end{document}